\documentclass[11pt,fleqn]{article}
\usepackage{amsfonts}
\newcommand{\ol}{\setlength{\itemsep}{0pt.}\begin{enumerate}}
\newcommand{\eol}{\end{enumerate}\setlength{\itemsep}{-\parsep}}
\newcommand{\ignore}[1]{}
\title{An upper bound for permanents of nonnegative matrices}
\author{Alex Samorodnitsky\thanks{School of Computer Science and
    Engineering, Hebrew University, Jerusalem, Israel.}}

\begin{document}
\date{}
\maketitle
 
 
\newtheorem{THEOREM}{Theorem}[section]
\newenvironment{theorem}{\begin{THEOREM} \hspace{-.85em} {\bf :} 
}%
                        {\end{THEOREM}}
\newtheorem{LEMMA}[THEOREM]{Lemma}
\newenvironment{lemma}{\begin{LEMMA} \hspace{-.85em} {\bf :} }%
                      {\end{LEMMA}}
\newtheorem{COROLLARY}[THEOREM]{Corollary}
\newenvironment{corollary}{\begin{COROLLARY} \hspace{-.85em} {\bf 
:} }%
                          {\end{COROLLARY}}
\newtheorem{PROPOSITION}[THEOREM]{Proposition}
\newenvironment{proposition}{\begin{PROPOSITION} \hspace{-.85em} 
{\bf :} }%
                            {\end{PROPOSITION}}
\newtheorem{DEFINITION}[THEOREM]{Definition}
\newenvironment{definition}{\begin{DEFINITION} \hspace{-.85em} {\bf 
:} \rm}%
                            {\end{DEFINITION}}
\newtheorem{EXAMPLE}[THEOREM]{Example}
\newenvironment{example}{\begin{EXAMPLE} \hspace{-.85em} {\bf :} 
\rm}%
                            {\end{EXAMPLE}}
\newtheorem{CONJECTURE}[THEOREM]{Conjecture}
\newenvironment{conjecture}{\begin{CONJECTURE} \hspace{-.85em} 
{\bf :} \rm}%
                            {\end{CONJECTURE}}
\newtheorem{MAINCONJECTURE}[THEOREM]{Main Conjecture}
\newenvironment{mainconjecture}{\begin{MAINCONJECTURE} \hspace{-.85em} 
{\bf :} \rm}%
                            {\end{MAINCONJECTURE}}
\newtheorem{PROBLEM}[THEOREM]{Problem}
\newenvironment{problem}{\begin{PROBLEM} \hspace{-.85em} {\bf :} 
\rm}%
                            {\end{PROBLEM}}
\newtheorem{QUESTION}[THEOREM]{Question}
\newenvironment{question}{\begin{QUESTION} \hspace{-.85em} {\bf :} 
\rm}%
                            {\end{QUESTION}}
\newtheorem{REMARK}[THEOREM]{Remark}
\newenvironment{remark}{\begin{REMARK} \hspace{-.85em} {\bf :} 
\rm}%
                            {\end{REMARK}}
 
\newcommand{\thm}{\begin{theorem}}
\newcommand{\lem}{\begin{lemma}}
\newcommand{\pro}{\begin{proposition}}
\newcommand{\dfn}{\begin{definition}}
\newcommand{\rem}{\begin{remark}}
\newcommand{\xam}{\begin{example}}
\newcommand{\cnj}{\begin{conjecture}}
\newcommand{\mcnj}{\begin{mainconjecture}}
\newcommand{\prb}{\begin{problem}}
\newcommand{\que}{\begin{question}}
\newcommand{\cor}{\begin{corollary}}
\newcommand{\prf}{\noindent{\bf Proof:} }
\newcommand{\ethm}{\end{theorem}}
\newcommand{\elem}{\end{lemma}}
\newcommand{\epro}{\end{proposition}}
\newcommand{\edfn}{\bbox\end{definition}}
\newcommand{\erem}{\bbox\end{remark}}
\newcommand{\exam}{\bbox\end{example}}
\newcommand{\ecnj}{\bbox\end{conjecture}}
\newcommand{\emcnj}{\bbox\end{mainconjecture}}
\newcommand{\eprb}{\bbox\end{problem}}
\newcommand{\eque}{\bbox\end{question}}
\newcommand{\ecor}{\end{corollary}}
\newcommand{\eprf}{\bbox}
\newcommand{\beqn}{\begin{equation}}
\newcommand{\eeqn}{\end{equation}}
\newcommand{\wbox}{\mbox{$\sqcap$\llap{$\sqcup$}}}
\newcommand{\bbox}{\vrule height7pt width4pt depth1pt}
\newcommand{\qed}{\bbox}
\def\sup{^}

\def\H{\{0,1\}^n}

\def\S{S(n,w)}

\def \E{\mathbb E}
\def \R{\mathbb R}
\def \Z{\mathbb Z}

\def\<{\left<}
\def\>{\right>}
\def \({\left(}
\def \){\right)}
\def \e{\epsilon}
\def \l{\lfloor}
\def \r{\rfloor}
\def \F{{\cal F}}
\def \G{{\cal G}}

\def \L{{\bigtriangledown}}

\def \a_i{\E d^2_{i,0}}
\def \b_i{\E d^2_{i,1}}

\def \O{\Omega(n,p)}
\def \U{U(n,p)}
\def \J{J(n,p)}
\def \l{\lambda}
\def \bl{{\bar \lambda}}
\def \by{{\bar y}}

\def \grad{\bigtriangledown}

\newcommand{\rarrow}{\rightarrow}
\newcommand{\lrarrow}{\leftrightarrow}

\begin{abstract}

\end{abstract}
A recent conjecture of Caputo, Carlen, Lieb, and Loss, and, independently,
of the author, states that the maximum of the permanent of a matrix whose
rows are unit vectors in $l_p$ is attained either for the identity matrix $I$
or for a constant multiple of the all-$1$ matrix $J$.

The conjecture is known to be true for $p = 1$ ($I$) and for $p \ge 2$ ($J$).

We prove the conjecture for a subinterval of $(1,2)$, and show the conjectured
upper bound to be true within a subexponential factor (in the dimension) for
all $1 < p < 2$. In fact, for $p$ bounded away from $1$, the conjectured upper
bound is true within a constant factor.

This leads to a mild (subexponential) improvement in
deterministic approximation factor for the permanent. We present an efficient
deterministic algorithm that approximates the permanent of a
nonnegative $n\times n$ matrix within $\exp\left\{n - O\(n/\log n\)\right\}$.

\section{Introduction}
Let $A = \(a_{ij}\)$ be an $n\times n$ matrix. The {\it permanent} of
$A$ is defined as
$$
per(A) = \sum_{\sigma \in S_n} \prod_{i=1}^n a_{i \sigma(i)}
$$
Here $S_n$ is the symmetric group on $n$ elements.

This paper investigates upper bounds on the permanent of matrices with
nonnegative entries. Bregman \cite{breg} resolved the Minc
conjecture and proved a tight upper bound on the permanent of a
zero-one matrix with given row sums. 
Here we are interested in upper bounds for matrices with general
nonnegative entries. (For related work see also \cite{Soules}
and the references there.)

More specifically, given $1 \le p \le \infty$,
we investigate the maximal possible value $U(n,p)$ of the permanent of
a matrix whose 
rows are unit vectors in $l^n_p$. We give an upper bound on $U(n,p)$ which
is tight up to a subexponential (in $n$) multiplicative factor. Since
the permanent is a multiliner function of its rows, this leads to
an upper bound on the permanent of an arbitrary real matrix, given the
$l_p$ length of its rows. 

Let us start with a conjecture claiming that there are
only two possible matrices on which the maximum of the permanent can
be attained. This conjecture is due to Caputo, Carlen, Lieb, and Loss
\cite{CLL}, and, independently, to the author.
\cnj
\label{cnj:l_p}
Let $1 \le p < \infty$. The maximum of the permanent of an $n\times n$
matrix whose rows are unit vectors in $l_p$ is attained in one of two
cases.
\begin{enumerate}
\item
On the identity matrix. In this case the
permanent is $1$.
\item
On a matrix all of whose entries are $n^{-1/p}$. In this case the
permanent is $\frac{n!}{n^{n/p}}$.
\end{enumerate}
In particular, the maximal possible value of the permanent is 
\beqn
\label{ineq:l_p}
U(n,p) = \mbox{max}\left\{1,
\frac{n!}{n^{n/p}}\right\}
\eeqn
\ecnj
Here are some preliminary remarks. Let the dimension $n$ be fixed. 
The function $f(p) = \frac{n!}{n^{n/p}}$
is increasing. Clearly $f(1) \le 1$ and $f(2) \ge 1$. It is
easy to compute the unique value of $p$, lying in $[1,2]$ for which
$f(p) = 1$, that is 
\beqn
\label{p_c}
p_c = \frac{n \log n}{\log n!}
\eeqn
Let $I$ denote the identity matrix, and $J$ the all-$1$ matrix. The
conjecture claims that $I$ is optimal for $p \in [1,p_c]$ and
$n^{-1/p} \cdot J$ is optimal for $p \in [p_c,\infty]$. 

In fact, it would suffice to prove the conjecture only for $p = p_c$. 
\lem
\label{lem:p_c}
\begin{itemize}
\item
Let $p_0 > 1$ be such that the matrix $I$ is optimal for
$p_0$. Then $I$ is the only optimal matrix for all $1 \le  p < p_0$.
\item
Let $p_0$ be such that the matrix $n^{1/p} \cdot J$ is optimal for
$p_0$. Then $n^{-1/p} \cdot J$ is the only optimal matrix for all $p > p_0$.
\end{itemize}
\elem
Let us now present the known results.
\begin{enumerate}
\item
The case $p = 1$ is trivial. For any $n$ only the identity matrix
is optimal, and $U(n,1) = 1$. 
\item
The conjecture is also known to be true for $p \ge 2$. In this case
the optimal matrix is $n^{-1/p} \cdot J$, and $U(n,p) =
\frac{n!}{n^{n/p}}$.  
Different proofs of this fact were given in \cite{SL98, Navon, SG01}. Later it
was pointed out \cite{G05} that this case was, essentially, already
dealt with in \cite{NN}. More specifically, the proof of \cite{SG01} is a
special case of an argument in \cite{NN} (Proposition 9.1.1, Appendix 1). 

To the best of our knowledge, the first published 
proof specifically treating this case appeared
recently in \cite{CLL}. Furthermore, this paper (independently) states
Conjecture~\ref{cnj:l_p}, attributing it also to P. Caputo.

Let us also mention that results in \cite{Fried78} imply
Conjecture~\ref{cnj:l_p} for $p \ge n$.
\item
The case $1 < p < 2$. This case seems to be the most
interesting. 

Clearly, one direction in (\ref{ineq:l_p}) is trivially
true: $U(n,p) \ge \mbox{max}\left\{1,\frac{n!}{n^{n/p}}\right\}$.
 
In the other direction, $U(n,p) \le U(n,1/2) = \frac{n!}{n^{n/2}}$.

This upper bound on $U(n,p)$ was improved in \cite{CLL}. They
show the function $U(n,p)$ to be logarithmically convex in $1/p$. This,
together with the known values $U(n,1) = 1$ and 
$U(n,2) = \frac{n!}{n^{n/2}}$, lead to an upper bound 
$$
U(n,p) \le \(\frac{n!}{n^{n/2}}\)^{2-2/p}
$$
\end{enumerate}

\noindent In this paper we show the conjecture to hold in the interval
$[1,p_0]$ where 
$$
p_0 = \frac{n\log n - (n-1)\log(n-1)}{\log n}
$$
For $n \ge 2$ holds $1 < p_0 < p_c \le 2$. 

It is interesting to
compare $p_0$ with $p_c$. We have $p_c \le \frac{\log
  n}{\log(n) - 1} = 1 + \frac{1}{\log(n) - 1}$. And $p_0 =
\frac{\log n + (n-1) \log \frac{n}{n-1}}{\log n} \ge \frac{\log n +
  (n-1)/n}{\log n} = 1 + \frac{1}{\log n} - \frac{1}{n \log n}$. 
Thus $p_c$ and $p_0$ are only about $\frac{1}{\log^2 n}$ apart. 

The proximity of $p_0$ and $p_c$, together with log-convexity of
$U(n,p)$, already suffice for giving an upper bound on $U(n,p)$ 
for all $p \in (1,2)$ which is tight up to a simply exponential factor
(in $n$). The approach we take will lead to a somewhat tighter
estimate, which has a subexponential error in the worst case.

Our main results are given in the following theorem.
\thm
\label{thm:main}
Let $n$ be fixed, and let $p_0 = \frac{n\log n - (n-1)\log(n-1)}{\log
  n}$.  
\begin{enumerate}
\item
The conjecture is true for $1 \le p \le p_0$. The identity
matrix is optimal for for $1 \le p \le
p_0$, and 
$$
U(n,p) = 1
$$
\item
For $p_0 \le p \le 2$ holds
$$
\max\left\{1,\frac{n!}{n^{n/p}}\right\} \le U(n,p) \le
  \exp\left\{(p-1)/p \cdot e^{1/(p-1)}\right\} \cdot \frac{n!}{n^{n/p}} 
$$
Observe, that this bound is $\exp\left\{n/\log n\right\}$-tight in the
worst case. 
For $p$ bounded away from $1$, this bound is tight within a constant factor.
\ignore{\item
The conjecture is true for $p\ge 2$. That is in this case holds
$$
U(n,p) = \frac{n!}{n^{n/p}}
$$
}
\end{enumerate}
\ethm
\ignore{
\rem
Our approach leads also to a new proof of the conjecture for the case
$p\ge 2$. We omit the details here (but see Remark~\ref{?}). 
\erem
}
\subsection{Approximating the permanent}
The original motivation for this study was computational.  
The goal is to construct an efficient deterministic algorithm that
approximates the permanent of a given nonnegative matrix within a reasonable
multiplicative factor. (A {\it randomized}
algorithm to approximate the permanent with arbitrary precision 
was constructed in \cite{JSV}.) 

In \cite{LSW} this problem was reduced to the case in which the input
matrix is doubly stochastic. This immeadiately gave an 
$\frac{n^n}{n!}$-approximation, since the permanent of a doubly
stochastic matrix lies between $\frac{n!}{n^n}$ and $1$. Here the
upper bound is trivial, while the lower bound is a deep theorem of
Egorychev \cite{Egor} and Falikman \cite{Falik}, proving a
conjecture by van der Waerden. In this light, it seems
natural to look fore more informative upper bounds, which could lead
to better approximation factors for the doubly-stochastic, and thus, for the
general case. 

Our results lead to an improvement of
$\exp\left\{O\(n/\log n\)\right\}$ in the approximation factor. We note that a
polynomial (in $n$) improvement in the approximation factor was
recently obtained in \cite{G2005}. 

The main tool is a permanental
inequality which might be of independent interest. This inequality is
an immediate consequence of Theorem~\ref{thm:main}.
\pro
Let $n \ge 2$ be an integer. Let $p_0 = \frac{n\log n - (n-1)\log(n-1)}{\log
n}$. Then for any stochastic $n \times n$ matrix $A = \(a_{ij}\)$ holds
$$
Per\(\(a^{1/p_0}_{ij}\)\) \le 1
$$
\epro
\cor
\label{cor:approx}
There is a deterministic polynomial-time algorithm to approximate the
permanent of a given nonnegative $n\times n$ matrix within a
multiplicative factor of $\frac{n^n}{n!} \cdot e^{-\Omega\(\frac{n}{log n}\)}$.
\ecor
\prf (Of the corollary)
It is sufficient to present an algorithm approximating the
permanent of a given doubly stochastic matrix within this factor. 

Let $q_0 = 1/p_0$. Assume $n \ge 5$. Let
$A$ be a doubly stochastic matrix. Let $\sigma \in S_n$ be a
permutation such that $\prod_{i=1}^n a_{i \sigma(i)}$ is
maximal. \footnote{Finding $\sigma$ amounts to finding a maximal
weight perfect matching in a given bipartite graph with $2n$
vertices, and can be done efficiently.}
Then there are two cases. 
\begin{itemize}
\item
$\prod_{i=1}^n a_{i \sigma(i)} \ge 2^{-n}$.
Then 
$$
2^{-n} \le \prod_{i=1}^n a_{i \sigma(i)} \le Per(A) \le 1
$$
\item
$\prod_{i=1}^n a_{i \sigma(i)} < 2^{-n}$. In this case, by the proposition,
$$
\frac{n!}{n^n} \le Per(A) \le 2^{(q_0-1)n} \cdot
Per\(\(a^{q_0}_{ij}\)\) \le  2^{(q_0-1)n} \le e^{-\Omega\(\frac{n}{log
    n}\)}
$$
\end{itemize}
\eprf

\subsection{Generalizations of Minc's conjecture to general
  nonnegative matrices}
The Minc conjecture, proved by Bregman, states that for a zero-one
matrix $A$ with $r_i$ ones in 
row $i$, $1\le i\le n$,
$$
per(A) \le \prod_{i=1}^n \(r_i !\)^{1/r_i}
$$
and equality holds if and only if $A$ is a block-diagonal matrix, and
all the blocks are all-$1$ matrices.\footnote{Up to a
  permutation of rows or columns.}
 
Let $\phi:~[0,1]
\rarrow [0,1]$ be a function taking $1/r$ to $1/\(r!\)^{1/r}$, for all integer
$r$. Given a matrix $A$ with entries in $[0,1]$, let $\phi(A)$ 
denote a matrix whose $(ij)$-th entry is $\phi\(a_{ij}\)$. 
Consider a stochastic matrix $A = \(a_{ij}\)$ whose $i$-th row
has entries with two possible values: $r_i$ entries with value $1/r_i$
and $n-r_i$ entries valued $0$. Then the Bregman bound implies
$$
per(\phi(A)) \le 1,
$$ 
and equality holds iff $A$ is a block-diagonal matrix with 
blocks which are constant multiples of all-$1$ matrices.

A natural way to extend $\phi$ to the whole interval $[0,1]$ is by
taking $\phi(x) = \Gamma\(1/x + 1\)^{-x}$, for all $0 < x \le 1$, and
setting $\phi(0) = 0$. The following conjecture generalizes the Minc
conjecture. 
\cnj
\label{cnj:minc_gen}
For any stochastic matrix $A$ holds
$$
per(\phi(A)) \le 1
$$
and equality holds iff $A$ is a block-diagonal matrix with
blocks which are constant multiples of all-$1$ matrices.
\ecnj
The function
$\phi = \Gamma\(1/x + 1\)^{-x}$ is strictly monotone and takes $[0,1]$ onto
$[0,1]$. It is also concave \cite{Moews}. 

Let $K = \left\{x\in \R^n;~\sum_{i=1}^n
\phi\(x_i\) \le 1\right\}$. This is a convex ball in $\R^n$
defining a norm $\|\cdot\|_K$. Consider the following optimization
problem: Choose $n$ unit vectors $x^{(1)}...x^{(n)}$ in $\R^n$ endowed
with the norm $\|\cdot\|_K$ as rows of a matrix so that the
permanent of this matrix is as large as
possible.~\footnote{Replacing permanent with determinant one arrives
to questions about the maximal volume subcube of $K$. These
questions are of interest in convex geometry \cite{Ball}. The two
contexts seem to be very different, however.}

An alternative way to state Conjecture~\ref{cnj:minc_gen} is to say
that all the optimal solutions to this optimization problem are
obtained as follows: partition $\{1...n\}$ into disjoint subsets
$S_1...S_k$. For each $j = 1...k$ choose all the vectors $x^{(i)}$,
$i\in S_j$, to be equal to $\frac{1}{|S_j|} \cdot {\bf 1}_{S_j}$, that is
be $\frac{1}{|S_j|}$ on the coordinates in $S_j$, and $0$ elsewhere.

The function $\phi$ and the norm it defines are somewhat compicated to
deal with. A natural ``easier'' family of norms to consider as a test
case are the $l_p$ norms, $1 \le p \le \infty$. This, in fact, was the
starting point of this study. 

We conclude the introduction by stating a conjecture which is a common
generalization of both Minc's conjecture and
Conjecture~\ref{cnj:l_p}. Following the discussion in
Lemma~\ref{lem:p_c}, Conjecture~\ref{cnj:l_p} is equivalent to $U(n,p_c) =
1$. Here $p_c = \frac{n \log n}{\log n!}$ is the `critical' value 
of $p$ for $n$-dimensional matrices. 

Let $p_c(r) = \frac{r \log r}{\log r!}$ for integer $r$. For $0 \le
r_1,r_2,...r_n \le n$ and $1 \le p_1,...,p_n < \infty$ let
$U(n;~r_1,...,r_n;~p_1,...,p_n)$ be the maximum of the permanent of an
$n\times n$ matrix whose $i$-th row is a unit vector in $l_{p_i}$
supported on at most $r_i$ non-zero coordinates. Then
\cnj
\label{cnj:common_gen}
$$
U(n;~r_1,...,r_n;~p_c(r_1),...,p_c(r_n)) \le 1
$$
\ecnj 
It is straightforward to check that for zero-one matrices this
conjecture is equivalent to the Minc conjecture. For $r_1=r_2=...=r_n
= n$ it reduces to Conjecture~\ref{cnj:l_p}.  
 
We remark that the proof of Theorem~\ref{thm:main} easily generalizes to
give 
$$
U(n;~r_1,...,r_n;~p_0(r_1),...,p_0(r_n)) \le 1
$$
where $p_0(r) = \frac{r\log r - (r-1)\log(r-1)}{\log
r}$. 

A word on our methods and an acknowledgement. Our proof of
Theorem~\ref{thm:main} proceeds along the lines of Bregman's proof of the Minc 
conjecture. A key inequality in that proof has to be replaced by a
more general inequality of \cite{baum}, quoted as  
Theorem~\ref{thm:baum} below. We are grateful to Leonid
Gurvits for directing us to this inequality.

\section{A recursive bound on U(n,p)}
Let $1 \le  p \le \infty$ be fixed. Let $q = 1/p$.

A vector $y = (y_1...y_n) \in \R^n$ is {\it
  stochastic} if its coordinates are nonnegative and sum to
$1$. Consider the following function defined on the set $\Delta$ of
stochastic vectors: 
$$
P(y) = \sum_{i=1}^n y^q_i \prod_{j\not = i} (1-y_j)^q
$$
This is a continuous bounded function which attains its maximum on
$\Delta$.
\dfn
$$
w(n,p) = \mbox{max}_{y\in \Delta} P(y)
$$
\edfn
The main claim of this section is:
\thm
\label{thm:perm_prod}
$$
\U \le \prod_{k=1}^n w(k,p)
$$
\ethm
\prf
The proof is by induction on $n$. For $n=1$, $U(1,p) = w(1,p) = 1$.
 
Consider an optimization problem 
$$
\mbox{Maximize} ~~~Per\(\lambda^q_{ij}\)
$$
Given  
$$
\lambda_{ij} \ge 0~~~~~~~~~~~~~\forall i~~\sum_{j=1}^n \lambda_{ij} =
1
$$
Clearly the optimal value here is $\U$. 

A key element of our proof is an inequality of \cite{baum}, which we
state next. 
\thm
\label{thm:baum}
Let $p(x,\lambda)$ be a nonnegative function defined on a space $X
\times \Lambda$ and let $\mu$ be a nonnegative weight function on
$X$. 

\noindent Let $P(\lambda) = \sum_{x \in X} \mu(x) p(x,\lambda)$, and
$Q(\lambda,{\bar \lambda}) = \sum_{x \in X} \mu(x) p(x,\lambda) \log
p(x,{\bar \lambda})$.

\noindent Then $Q(\lambda, {\bar \lambda}) \ge Q(\lambda, \lambda)$ implies
$P(\bl) > P(\l)$ unless $p(x,\l) = p(x,\bl)$ for all $x$ with $\mu(x)
> 0$.
\ethm 

Now we apply Theorem~\ref{thm:baum} in our setting.

Let $X = S_n$ be the symmetric group on $n$ elements, and
$\Lambda$ be the set of all stochastic matrices $(\l_{ij})$. Let
$\mu(\sigma) = 1$ for all permutations $\sigma \in S_n$ and let $p(\sigma,\l)
= \prod_{i=1}^n \l^q_{i,\sigma(i)}$, for $\sigma \in S_n$ and $\l \in
\Lambda$.  Then $P(\l) =  \sum_{x \in X} \mu(x) p(x,\lambda) =
Per\(\l^q_{ij}\)$.  

Let $\l[i,j]$ be the $(n-1) \times (n-1)$ matrix
obtained from $\l$ by deleting $i$-th row and $j$-th column. Let
$\l^q[i,j]$ be the matrix obtained from $\l[i,j]$ by raising each
entry to $q$-th power. Let $\bl
= (\bl_{ij})$ with 
$$
\bl_{ij} = \frac{\l^q_{ij} Per\(\l^q[i,j]\)}{Per\(\l^q_{ij}\)}
$$ 
The following lemma is a direct consequence of Theorem~\ref{thm:baum}.
\lem
\label{lem:baum}
$$
Per\(\bl^q_{ij}\) \ge Per\(\l^q_{ij}\)
$$
\elem
\prf
Consider the optimization problem of maximizing $Q\(\l,\bl\)$ given
$\l$. We have 
$$
Q(\l,\bl) = \sum_{\sigma \in S_n} p(\sigma, \lambda)
\log p(\sigma, \bl) = q \cdot \sum_{\sigma \in S_n} p(\sigma, \lambda)
\sum_{i=1}^n \log \bl_{i\sigma(i)} = q \cdot \sum_{i,j=1}^n
\l^q_{ij} Per\(\l^q[i,j]\) \log \bl_{ij}
$$
The constraints on $\bl$ are that it is a stochastic matrix. Therefore we
have $n$ independent optimization problems of the form:
$$
\mbox{Maximize}~~\sum w_j \log y_j~~~~~~~~~~\mbox{Given}~~y_j\ge 0,~~\sum y_j
= 1,
$$
where $w_j$ are nonnegative constants. Assuming not all $w_j$ are
zero, which we may and will do in our case, the only solution of this
problem is $y_j = \frac{w_j}{\sum_k w_k}$. This is a simple
consequence of the concavity of the logarithm.

Fixing $1 \le i\le n$, and substituting $w_j = \l^q_{ij}
Per\(\l^q[i,j]\)$ and $y_j = \bl_{ij}$, we see that optimal $\bl$ is
given by $\bl_{ij} = \frac{\l^q_{ij}
  Per\(\l^q[i,j]\)}{Per\(\l^q_{ij}\)}$. The claim of the lemma now
follows from Theorem~\ref{thm:baum}.
\eprf

Now, following \cite{breg}, we write
$$
Per\(\l^q_{ij}\) \le Per\(\bl^q_{ij}\) = \sum_{\sigma \in S_n} \prod_{i=1}^n
\bl^q_{i\sigma i} = 
$$
$$
\sum_{\sigma \in S_n} \prod_{i=1}^n
\(\frac{\l^q_{i\sigma(i)}
  Per\(\l^q[i,\sigma(i)]\)}{Per\(\l^q_{ij}\)}\)^q =
\frac{1}{Per^{qn}\(\l^q_{ij}\)} \cdot \sum_{\sigma \in S_n}
\prod_{i=1}^n \l^{q^2}_{i\sigma(i)} Per^q\(\l^q[i,\sigma(i)]\)
$$
 
Let $(\l^q_{ij})$ be an optimal matrix, that is
$Per\((\l^q_{ij})\) = \U$. Then 
$$
\U^{qn + 1} \le \sum_{\sigma \in S_n}
\prod_{i=1}^n \l^{q^2}_{i\sigma(i)} Per^q\(\l^q[i,\sigma(i)]\)
$$ 

Consider the matrix $\l[i,j]$. This is an $(n-1)\times (n-1)$
matrix with row sums $r_k = 1 - \l_{kj}$, for $k = 1...n$, $k \not =
i$. Let $R$ be 
the $(n-1)\times (n-1)$ diagonal matrix with $1/r_k$ on the
diagonal. Then $(a_{ij}) =  R \cdot \l[i,j]$ is a stochastic matrix,
and therefore, by induction hypothesis,
$Per(a^q_{ij}) \le U(n-1,p)$. 

This means $Per\(\l^q[i,j]\) \le
U(n-1,p) \cdot \prod_{k\not = i} (1-\l_{kj})^q$. Substituting this in
the inequality above, we obtain
$$
\U^{qn + 1} \le U(n-1,p)^{qn} \cdot \sum_{\sigma \in S_n}
\prod_{i=1}^n \(\l^{q^2}_{i\sigma(i)} \prod_{k\not = i}
(1-\l_{k\sigma(i)})^{q^2}\) = 
$$ 
$$
 U(n-1,p)^{qn} \cdot \sum_{\sigma \in S_n} \(\prod_{i=1}^n
 \l_{i\sigma(i)} \cdot \prod_{k,j:~\sigma(k)\not = j}
 (1-\l_{kj})\)^{q^2} =
$$
$$
U(n-1,p)^{qn} \cdot \prod_{i,j} \(1-\l_{ij}\)^{q^2} \cdot
\sum_{\sigma \in S_n} \prod_{i=1}^n
\(\frac{\l_{i\sigma(i)}^q}{\(1-\l_{i\sigma(i)}\)^q}\)^q
$$ 
The third term in this expression is the permanent of a matrix
$\(a^q_{ij}\)$, where $a_{ij} = \frac{\l_{ij}^q}{\(1-\l_{ij}\)^q}$.

Let $r_i = \sum_{j=1}^n a_{ij} = \sum_{j=1}^n
\frac{\l_{ij}^q}{\(1-\l_{ij}\)^q}$ be the row sums of this
matrix. Then, $Per(a^q_{ij}) \le \U \cdot 
\prod_{i=1}^n r^q_i$. Substituting in the inequality above gives
$$
\U^{qn} \le U(n-1,p)^{qn} \cdot \prod_{i,j} \(1-\l_{ij}\)^{q^2} \cdot 
\(\prod_{i=1}^n \sum_{j=1}^n \frac{\l_{ij}^q}{\(1-\l_{ij}\)^q}\)^q
$$  
Taking $q$-th roots of both sides this simplifies to
$$
\U^n \le U(n-1,p)^n \cdot \prod_{i,j} \(1-\l_{ij}\)^q \cdot 
\prod_{i=1}^n \sum_{j=1}^n \frac{\l_{ij}^q}{\(1-\l_{ij}\)^q}
$$
Let $\l_i$ be the $i$-th row vector of $\l$. Since $\l$ is a
stochastic matrix, $\l_i$ is a stochastic vector. We have 
$$
\prod_{i,j} \(1-\l_{ij}\)^q \cdot 
\prod_{i=1}^n \sum_{j=1}^n \frac{\l_{ij}^q}{\(1-\l_{ij}\)^q} =
\prod_{i=1}^n P(\l_i) \le w(n,p)^n
$$
Therefore $U(n,p) \le U(n-1,p) \cdot w(n,p)$. The claim now follows
from the induction hypothesis
$$
U(n,p) \le U(n-1,p) \cdot w(n,p) \le  w(n,p) \cdot \prod_{k=1}^{n-1} w(k,p) 
  = \prod_{k=1}^n w(k,p)
$$
\eprf
\section{Proofs of the main results}
Our first order of business is to determine $w(k,p)$, for $1\le k\le
n$. Let $1 < p < 2$ be fixed, and let $q = 1/p$.

Let $\theta(k) = k \cdot \(\frac{(k-1)^{k-1}}{k^k}\)^q$ for integer $k \ge 2$
and let $\theta(1) = 1$. 
\thm
\label{thm:wkp}
Fix $k \ge 1$. The maximum of $P(y) = \sum_{i=1}^k y^q_i \prod_{j\not =
  i} (1-y_j)^q$ is 
attained either at a standard basis vector and then $w(k,p) =
\theta(1) = 1$, or at the all-$1/k$ vector, in which case $w(k,p) = \theta(k)$. 
\ethm
\ignore{
\rem
This result could be interpreted as a statement about permanents. For
a stochastic vector $y$, let $\l = \l(y)$ be a matrix whose first row
is given by $\l_{1j} = y_j$, and all the other rows are equal and are
given by $\l{ij} = (1-y_j)/(k-1)$, for $i>1$. Then $\l$ is a doubly
stochastic matrix, and $Per\(\l^q\) = \frac{(k-1)!}{(k-1)^{q(k-1)}} \cdot
P(y)$.

Going through the stochastic vectors $y \in \R^k$ leads to a
restricted class of doubly stochastic matrices 
containing two special matrices. The first matrix is $A = \(1/k\)
\cdot J$. The second matrix $B$
is a block-diagonal matrix, with two blocks of size $1$ and
$k-1$ respectively. Both blocks are constant multiples of an all-$1$ matrix.
The theorem says that for each $q$ the  
maximum of the function $Per\(\l^q\)$ on this class of matrices is
attained either at $A$ or at $B$.
\erem
}
The proof of Theorem~\ref{thm:wkp} is technical and is relegated to
Appendix.

We briefly discuss the claim of the theorem. Let $I_k$ be the $k\times
k$ identity matrix. Let $J_k$ denote the
matrix $k^{-1/p} \cdot J$, where $J$ is the all-$1$ $k\times k$
matrix. Note that $\theta(k) =
\frac{per\(J_k\)}{per\(J_{k-1}\)}$. Therefore the theorem, combined
with Theorem~\ref{thm:perm_prod}, says that for any $k \ge 2$ 
$$
\frac{U(k,p)}{U(k-1,p)} \le
\max\left\{\frac{per\(I_k\)}{per\(I_{k-1}\)},\frac{per\(J_k\)}{per\(J_{k-1}\)}\right\}       
$$ 
Let us observe that this inequality agrees well with Conjecture~\ref{cnj:l_p}.

The last step before the proof of Theorem~\ref{thm:main} is
Lemma~\ref{lem:p_c}, which we prove now.

\noindent \prf (Lemma~\ref{lem:p_c}) 

The following notation will be convenient. For $1 \le
p \le \infty$, let $\O$ be
the set of $n \times n$ matrices whose rows are unit vectors in $l_p$.

We need a following well-known fact. Let $1 \le p < p' \le
\infty$. Let $a$ be a vector in $\R^n$. Then
\beqn
\label{norms}
1 \le \frac{\|a\|_p}{\|a\|_{p'}} \le n^{\frac 1p - \frac{1}{p'}}
\eeqn
Equality on the left is possible only for a multiple of a standard basis
vector, and equality on the right is possible only for a multiple of
the all-$1$ vector. 

Let $p_0$ be such that the matrix $I$ is optimal for
$p_0$. Let $p < p_0$. Let $A \in \O$ with rows $a_1...a_n$. 
 Let $D = \(d_{ij}\)$ be a
diagonal matrix with $d_{ii} = \frac{\|a_i\|_p}{\|a_i\|_{p_0}}$. Then $DA$
is in $\Omega(n,p_0)$ and therefore
$$
per(A) = per\(D^{-1} \cdot (DA)\) = per\(D^{-1}\) \cdot per(DA) = 
\prod_{i=1}^n \frac{\|a_i\|_{p_0}}{\|a_i\|_p} \cdot per(DA)
$$
$$
\le \prod_{i=1}^n \frac{\|a_i\|_{p_0}}{\|a_i\|_p} \cdot per(I) =
\prod_{i=1}^n \frac{\|a_i\|_{p_0}}{\|a_i\|_p} \le 1
$$
By (\ref{norms}) 
equality is only possible if all the rows $a_i$ are standard basis
vectors, and $A$ is the identity matrix, up to permuting coordinates.

This proves the first claim of the lemma. The proof of the second claim
proceeds along similar lines, using second half of
inequality (\ref{norms}). We omit the details.
\eprf

\noindent \prf (Theorem~\ref{thm:main})

Fix $p = p_0 = \frac{n\log n - (n-1)\log(n-1)}{\log n}$. Let $q =
1/p$. The value of $p$ is chosen precisely so that $\theta(n) = n \cdot
\(\frac{(n-1)^{n-1}}{n^n}\)^q = 1$.

By Theorem~\ref{thm:perm_prod}, Theorem~\ref{thm:wkp}, and
Lemma~\ref{lem:mon_theta}   
$$
U(n,p) \le \prod_{k=1}^n w(k,p) \le \prod_{k=1}^n
\max\left\{1,\theta(k)\right\} \le 
\(\max\left\{1,\theta(n)\right\}\)^n = 1
$$
Therefore $I$ is optimal for $p = p_0$. Lemma~\ref{lem:p_c} completes
the proof of the first claim of the theorem.

Now, to the second claim. Fix $p \in (1,2)$. Let $q = 1/p$. By
Lemma~\ref{lem:mon_theta} there is an integer $k_0$ such that
$\theta(k) < 1$ for $k \le k_0$ and $\theta(k) \ge 1$ for $k >
k_0$. Since $\theta(k) = \frac{per\(J_k\)}{per\(J_{k-1}\)}$, this
means that $per\(J_{k_0}\) = \prod_{k=1}^{k_0} \theta(k) = \min_{k\ge 1}
per\(J_k\)$.

Therefore,
$$
U(n,p_0) \le \prod_{k=1}^n w(k,p_0) \le \prod_{k=1}^n
\max\left\{1,\theta(k)\right\} = \prod_{k=2}^n
\max\left\{1,\frac{per\(J_k\)}{per\(J_{k-1}\)}\right\} = 
$$
$$
\frac{per\(J_n\)}{per\(J_{k_0}\)} = \frac{per\(J_n\)}{\min_{k\ge 1}
  per\(J_k\)}
$$
It remains to estimate the denominator on the right.

We have
$$
\min_{k\ge 1} per\(J_k\) = \min_{k\ge 1} \frac{k!}{k^{qk}} \ge
\min_{k\ge 1} \frac{k^{(1-q)k}}{e^k} \ge \min_{x\ge 1} \frac{x^{(1-q)x}}{e^x}
$$
where in the last inequality an integer variable
$k$ is replaced with a real variable $x$. A simple analysis gives that
the minumum on the right hand side is attained for $x =
\exp\left\{q/(1-q)\right\} = \exp\left\{1/(p-1) \right\}$ and equals  
$\exp\left\{-(p-1)/p \cdot e^{1/(p-1)}\right\}$. 

Therefore
$$
U(n,p) \le \exp\left\{(p-1)/p \cdot e^{1/(p-1)}\right\} \cdot
per\(J_n\) = \exp\left\{(p-1)/p \cdot e^{1/(p-1)}\right\} \cdot
\frac{n!}{n^{n/p}} 
$$  
This completes the proof of the second claim and of the theorem.
\eprf
\ignore{
\rem
The bound in theorem~\ref{perm_prod}
can easily be replaced with  
$$
\U \le U(m,p) \cdot \prod_{k=m+1}^n w(k,p),
$$
for any $1 \le m\le n$. This implies, for instance, that if we prove
the conjecture for $n = 3$ this would immediately give optimality of
$\J$ for any $n$ and for $p$ greater than the critical $p$ for $n=3$,
and so on.
\erem
}

\section{Appendix: A Proof of Theorem~\ref{thm:wkp}}
We start with a useful property of the function $\theta$. Let $1/2 < q
< 1$ be a real number.
\lem
\label{lem:mon_theta}
Let $k \ge 1$ and consider the continuous function $\theta(x) = x \cdot
\(\frac{(x-1)^{x-1}}{x^x}\)^q$ of a real variable $x$ on the interval
$[1,k]$. If $x_0$ is a point of maximum of $\theta$ then $x_0 = 1$ or
$x_0 = k$.
\elem
\prf
It is convenient to deal with $f(x) = \ln(\theta(x)) = \ln x - q \cdot
\(x\ln x - (x-1) \ln (x-1) \)$. The derivative $f'(x) =
\frac{1}{x} - q \ln \frac{x}{x-1} = \frac{q}{x} \cdot \(\frac{1}{q} -
x \ln \frac{x}{x-1}\)$. 

Consider the function $g(x) = x \ln \frac{x}{x-1}$ on
$[1,\infty)$. The derivative 
$g'(x) = \ln \frac{x}{x-1} - \frac{1}{x-1} = \ln \(1 + \frac{1}{x-1}\)
- \frac{1}{x-1}$ is strictly negative. At the endpoints, $g(1) =
\infty$ and $g(\infty) = 1$. Therefore on $[1,\infty)$ the function
$g$ decreases from $\infty$ to $1$. Since $\frac{1}{q} = p >
  1$ this means that there exists a
positive real number $x_q > 1$ depending only on $q$ such that $f' <
0$ for $1 \le x < x_q$, $f'\(x_q\) = 0$, and $f'(x) > 0$ for $x >
x_q$. 

Consequently, $f$ is unimodal on $[1,\infty)$ with minimum in
  $x_q$. The claim of the lemma follows.
\eprf

The proof of the theorem proceeds by induction on $k$. For $k=1$ the
claim holds trivially. For $k = 2$ we have
$$
P(y) = P(y_1,1-y_1) = y^{2q}_1 + (1-y_1)^{2q}
$$
For $q > 1/2$, the function $f(x) = x^{2q} + (1-x)^{2q}$ attains its
maximum on $[0,1]$ at $0$ and at $1$. This means that the
points of maximum of $P$ are standard basis vectors, and the claim holds.

Assume the theorem is true for $2\le l < k$.

Let $y^* \in \Delta$ be a point at which $P$ attains maximum. If $y^*$ has
$1 < l < k$ non-zero coordinates, then the induction hypothesis implies
$y^*$ is the all-$1/l$ vector. This is to say $P\(y^*\) =
\theta(l)$. However, Lemma~\ref{lem:mon_theta} showed $\theta(l) <
\max\left\{1,\theta(k)\right\}$, reaching a contradiction.

Therefore either $y^*$ is a standard basis vector, in which case we
are done, or 
$y^*$ is an interior point of $\Delta$. This is the remaining case. We
will assume that $y^*$ is not the all-$1/k$ vector and reach a contradiction.

Since $y^*$ is an interior extremum point, we can use the first and the second
order optimality conditions on the gradient and the Hessian of $P$ at
$y^*$ to obtain information about $y^*$.

Let $s_i(y) = y^q \prod_{j\not = i} \(1-y_j\)^q$, for $i = 1...k$. Of
course $P = \sum_{i=1}^k s_i$.
\lem
\label{lem:grad}
For all $i = 1...k$ 
$$
s_i\(y^*\) = y^*_i P\(y^*\)
$$
\elem
\prf
We have $\frac{\partial s_i}{\partial y_i} = \frac{q s_i}{y_i}$ and,
for $j \not = i$, $\frac{\partial s_i}{\partial y_j} = -\frac{q
  s_i}{1 - y_j}$. Therefore
$$
\frac{\partial P}{\partial y_j} = \sum_{i=1}^n \frac{\partial
  s_i}{\partial y_j} = \frac{\partial s_j}{\partial y_j} +
\frac{\partial P - s_j}{\partial y_j} = q\cdot \(\frac{s_j}{y_j} -
\frac{P - s_j}{1-y_j}\) = q \cdot \frac{s_j - y_j P}{y_j(1-y_j)}
$$
The first order optimality conditions for $y^*$ say that there is a
constant $\lambda$ such that for all $j=1...k$ holds $\frac{\partial
  P}{\partial y_j}\(y^*\) = \lambda$. This means that for $j=1...k$
holds $s_j\(y^*\) - y^*_j P\(y^*\) = \frac{\lambda}{q}y^*_j\(1-y^*_j\)$. 

Summing over $j$ we obtain 
$$
\frac{\lambda}{q} \cdot \sum_{j=1}^k y^*_j\(1-y^*_j\) = 0, 
$$
implying $\lambda = 0$. That is, for all $j = 1...k$ holds $s_j\(y^*\)
= y^*_j P\(y^*\)$. 
\eprf
\cor
\label{cor:ab}
The coordinates of $y^*$ have two distinct values $a$ and $b$ with $a
< 1- q < b$.
\ecor
\prf
Let $i \not = j$ be two distinct indices. By the lemma at $y^*$ we
have $s_i = y^*_i P$ and $s_j = y^*_j P$. This implies 
$$
\frac{y^*_i}{y^*_j} = \frac{s_i}{s_j} = \frac{\(y^*_i\)^q
  \(1-y^*_j\)^q}{\(y^*_j\)^q \(1-y^*_i\)^q}
$$ 
This means $\(y^*_i\)^{1-q}\(1-y^*_i\)^q =
\(y^*_j\)^{1-q}\(1-y^*_j\)^q$. Let $f(x) = x^{1-q}(1-x)^q$. We have
shown that $f\(y^*_i\) = f\(y^*_j\)$. Since the argument does not
depend on the choice of $i$ and $j$, this implies $f$ has the same
value on all $y^*_i$, $i = 1...k$. 

The function $f$ is
a concave function on $[0,1]$ vanishing at the endpoints, with maximum
at $1-q$. Therefore $f$ takes each value at most twice, at two points
lying on different sides of $1-q$. Bearing in mind that $y^*$ is not a
constant vector, the claim of the corollary follows.
\eprf

Next, we compute the Hessian of $P$. We have, for $i \not = j \not = t$
$$
\frac{\partial^2 s_i}{\partial y^2_i} = -\frac{q(1-q)
  s_i}{y^2_i};~~~~~\frac{\partial^2 s_i}{\partial y^2_j} = -\frac{q(1-q)
  s_i}{\(1-y_j\)^2};
$$
$$
\frac{\partial^2 s_i}{\partial y_i \partial
  y_j} = -\frac{q^2 s_i}{y_i\(1-y_j\)};~~~~~\frac{\partial^2
  s_i}{\partial y_j \partial y_t} = \frac{q^2 s_i}{\(1-y_j\)\(1-y_t\)}
$$ 
Let $H = H(y)$ be the Hessian of $P$ at $y$. Then 
$$
H(j,j) = \frac{\partial^2 P}{\partial y^2_j} = \sum_{i=1}^k
\frac{\partial^2 s_i}{\partial y^2_j} = -q(1-q) \cdot
\(\frac{s_j}{y^2_j} + \frac{P-s_j}{(1-y_j)^2}\) 
$$   
Similarly
$$
H(j,t) = \frac{\partial^2 P}{\partial y_j \partial y_t} = \sum_{i=1}^k
\frac{\partial^2 s_i}{\partial y_j \partial y_t} = \frac{\partial^2
  s_j}{\partial y_j \partial y_t} + \frac{\partial^2
  s_t}{\partial y_j \partial y_t} + \frac{\partial^2
  \(P - s_j - s_t\)}{\partial y_j \partial y_t} = 
$$
$$
-q^2\cdot
\(\frac{s_j}{y_j (1-y_t)} + \frac{s_t}{y_t (1-y_j)} \) + q^2 \cdot
\frac{P - s_j - s_t}{(1-y_j)(1-y_t)} 
$$
At $y^*$ we have $s_i = y^*_i P$ for all $i = 1...k$. Therefore for $H
= H\(y^*\)$ we have
$$
H(j,j) = -q(1-q) \cdot \frac{P}{y^*_j(1-y^*_j)}
$$
and
$$
H(j,t) = q^2 P \cdot \(\frac{1}{(1-y^*_j)(1-y^*_t)} - \frac{1}{(1-y^*_j)} -
\frac{1}{(1-y^*_t)}\) = -q^2 \cdot \frac{P}{(1-y^*_j)(1-y^*_t)}
$$
\lem
\label{lem:hess}
$y^*$ has only one coordinate with value $b$. (And therefore $k-1$
coordinates with value $a$.)
\elem
\prf
We can write the Hessian at $y^*$ as $H = -qP\cdot \(A + D\)$, where
$A$ is a rank-$1$ matrix with $a_{ij} =
\frac{q}{\(1-y^*_i\)\(1-y^*_j\)}$, and $D$ is a diagonal matrix with
$d_{ii} = \frac{1-q-y^*_i}{y^*_i\(1-y^*_i\)^2}$.

The second order optimality conditions for $y^*$ say that $H$ is 
negative semidefinite on the subspace $V$ of the vectors in $\R^k$
orthogonal to the all-$1$ vector. This means that the matrix $B = A +
D$ is positive semidefinite on $V$.

Assume for the moment that $y^*$ has two $b$-valued coordinates. Let
these be the first two coordinates. This means that 
$\left[\begin{array}{cc}
a_{11} & a_{12} \\
a_{21} & a_{22} 
\end{array}
\right] = \left[\begin{array}{cc}
1/(1-b)^2 & 1/(1-b)^2 \\
1/(1-b)^2 & 1/(1-b)^2 
\end{array}
\right]
$,
and $\left[\begin{array}{cc}
d_{11} & d_{12} \\
d_{21} & d_{22} 
\end{array}
\right] = \left[\begin{array}{cc}
\frac{1-q-b}{b(1-b)^2} & 0 \\
0 & \frac{1-q-b}{b(1-b)^2} 
\end{array}
\right]
$.
Note, that since $b > 1-q$, the diagonal values of the second matrix
are negative.

Now, let $v \in V$, $v = (1,-1,0,\ldots,0)$. Then clearly 
$$
v B v^t = 2 \frac{1-q-b}{b(1-b)^2} < 0,
$$
contradicting positive semidefinitness of $B$. This means that $y^*$
has only one coordinate valued $b$. 
\eprf

Consider the set $\Delta_1 \subset \Delta$ 
of stochastic vectors $y$ with $y_2 = ... = y_k =
\frac{1-y_1}{k-1}$. The preceding lemma implies that there is a maximum point
$y^*$ of $P$ in $\Delta_1$. Moreover $b = y^*_1 > 1/k$. 

$P$, restricted to $\Delta_1$, is a function of one variable $x = y_1$
and is given by
$$
P(x) = x^q \(1-\frac{1-x}{k-1}\)^{(k-1)q} + (k-1)\(\frac{1-x}{k-1}\)^q
  (1-x)^q \(1-\frac{1-x}{k-1}\)^{(k-2)q} = 
$$
$$
\frac{1}{(k-1)^{(k-1)q}}
      \cdot \(x^q (k-2+x)^{(k-1)q} + (k-1)(1-x)^{2q}(k-2+x)^{(k-2)q}\)
$$
We will show that on the interval $\left[1/k,1\right]$ 
this function attains its maximum either at $1/k$ or at
$1$. This means, recalling $y^*_1 > 1/k$, that $y^*$ is a standard basis 
vector. This is a contradiction to previous assumptions, and will
complete the proof of the theorem. 
\lem
\label{lem:one_var}
Let $k \ge 3$ be an integer, let $1/2 < q < 1$ be a real number, and
let $f$ be a function on $[1/k,1]$ given by  
$$
f(x) = x^q (k-2+x)^{(k-1)q} + (k-1)(1-x)^{2q}(k-2+x)^{(k-2)q}
$$
Then $f$ attains its maximum either at $1/k$ or at
$1$.
\elem
\prf
We compute the derivative of $f$.
$$
f'(x) = qx^{q-1} (k-2+x)^{(k-1)q}~+~(k-1)q x^q (k-2+x)^{(k-1)q - 1}~-
$$
$$
2(k-1)q (1-x)^{2q-1} (k-2+x)^{(k-2)q}~+~(k-1)(k-2)q (1-x)^{2q}
(k-2+x)^{(k-2)q -1} = 
$$
$$
q(k-2+kx)(k-2+x)^{(k-2)q-1}x^{q-1}~\cdot 
~\((k-2+x)^q~-~(k-1)x^{1-q}(1-x)^{2q-1}\) 
$$  

This means that the sign of $f$ is determined by the sign of
$(k-2+x)^q~-~(k-1)x^{1-q}(1-x)^{2q-1}$. 

Since $t \mapsto t^q$ is monotone
increasing, we can, as well, check the sign of 
$$
h(x) = (k-2+x)~-~(k-1)^{1/q} x^{1/q-1} (1-x)^{2-1/q}
$$
The function $h(x)$ is strictly convex on $[1/k,1]$, with $h\(1/k\) = 0$ and
$h(1) = k-1 > 0$.  

Therefore, there are two possible options.
\begin{itemize}
\item
$h > 0$ on $(1/k,1]$. This means that $f$ attains its maximum at $1$.
\item
There is a point $x \in \(1/k, 1\)$ such that 
$h < 0$ on $(1/k,x)$ and $h > 0$ on 
  $(x,1)$. This means that $f$ attains its maximum at one of the
  endpoints $1/k$ or $1$, and we are done. 
\end{itemize}  
\eprf

This completes the proof of Theorem~\ref{thm:wkp}.

\section{Acknowledgements}
We are grateful to Leonid Gurvits for several very helpful
discussions. We also thank Shmuel Friedland, Nati Linial, and Michael
Navon for valuable conversations.


\begin{thebibliography}{99}

\bibitem{baum}
L. E. Baum, T. Petrie, G. Soules, and N. Weiss, {\it A maximization
technique occurring in the statistical analysis of probabilistic
functions of Markov chains}, Ann. of Math. Stat., vol. 41, 1, 1970
pp. 164-171.

\bibitem{Ball}
K. Ball, {\it An elementary introduction to modern convex geometry}, Flavors
of Geometry, MSRI Publications, Vol.31, 1997.

\bibitem{breg}
L. M. Bregman, {\it Some properties of nonnegative matrices and their
  permanents}, Soviet Math. Dokl., vol. 14, 4, 1973, pp. 945-949. 

\bibitem{CLL}
E. Carlen, M. Loss, and E. H. Lieb, {\it An inequality of Hadamard
  type for permanents}, Mathematics ArXiv NT/0508096, 2005.

\bibitem{Egor} 
G.P. Egorychev, {\it The solution of van der Waerden's problem 
for permanents}, Advances in Math., 42, 299-305, 1981.

\bibitem{Falik}
D. I. Falikman, {\it Proof of the van der Waerden's conjecture on the
permanent of a doubly stochastic matrix}, Mat. Zametki
29, 6: 931-938, 957, 1981, (in Russian).

\bibitem{Fried78}
S. Friedland, {\it A study of the van der Waerden conjecture and its
  generalizations}, Linear and Multilinear Algebra, 6, 1978, pp. 123-143.

\bibitem{G2005}
L. Gurvits, {\it Hyperbolic polynomials approach to van der
  Waerden/Schrijver-Valiant like conjectures: sharper bounds, simpler
  proofs and algorithmic applications}, to appear in STOC 2006.

\bibitem{G05}
L. Gurvits, personal communication.

\bibitem{SG01}
L. Gurvits, A. Samorodnitsky, an unpublished manuscript, 2001.

\bibitem{JSV} M. Jerrum, A. Sinclair and E. Vigoda, {\it A polynomial-time
approximation algorithm for the permanent of a matrix with
non-negative entries}, Proc. 33 ACM Symp. on Theory of Computing,
ACM, 2001.

\bibitem{SL98}
N. Linial, A. Samorodnitsky, an unpublished manuscript, 1998.

\bibitem{LSW}
N. Linial, A. Samorodnitsky and A. Wigderson, {\it A deterministic strongly 
polynomial algorithm for matrix scaling and approximate permanents},
Combinatorica, vol. 20, 4, 2000. 

\bibitem{Moews}
D. Moews, {\it $\Gamma(x+1)^{1/x}$ is concave}, personal communication
to the author of \cite{Soules}.

\bibitem{Navon}
M. Navon, {\it Some notes on the permanent}, a project report,
submitted as a part of MSc thesis requirements, Hebrew University, 2004.

\bibitem{NN}
Y. Nesterov, A. Nemirovski, {\bf Interior Point Polynomial Methods in
  Convex Programming}, SIAM, Philadelphia, 1994.

\bibitem{Soules}
G. W. Soules, {\it New permanental upper bounds for nonnegative
  matrices}, Linear and Multilinear Algebra 51, 2003, pp. 319-337.


\end{thebibliography}
\end{document}